\documentclass[a4paper]{article}
\usepackage{amsmath}
\usepackage{amssymb}
\usepackage{theorem}
\usepackage{latexsym}
\usepackage{subfig}
\usepackage{subfig}
\usepackage{graphicx,color}
\usepackage{multirow}
\usepackage{colortbl}
\usepackage{rotating}
\usepackage[table]{xcolor}
\usepackage{hhline}
\usepackage{floatrow}
\usepackage{mathrsfs}

\DeclareMathAlphabet{\mathantt}{OT1}{antt}{li}{it}
\DeclareMathAlphabet{\mathpzc}{OT1}{pzc}{m}{it}

\definecolor{lightgray}{gray}{0.95}
\numberwithin{equation}{section}

\newfloatcommand{capbtabbox}{table}[][\FBwidth]

\newtheorem{defi}{Definition}

\input amssym.def
\input amssym.tex

\def\qed{\raise1pt\hbox{\vrule height5pt width5pt depth0pt}}
\title{{\bf Basins of attraction in forced systems\\with time-varying dissipation}}
\author{James A. Wright{\small$^1$}, Jonathan H.B. Deane{\small$^1$}, Michele Bartuccelli{\small$^1$}, Guido Gentile{\small$^2$} \\ \\
{\small $^1$Department of Mathematics, University of Surrey, Guildford, GU2 7XH, UK} \\
{\small $^2$Dipartimento di Matematica e Fisica, Universit\`{a} di Roma Tre, 00146 Roma, Italy} \\
{\small Email: j.wright, j.deane, m.bartuccelli@surrey.ac.uk, gentile@mat.uniroma3.it}}
\date{}
\begin{document}
\maketitle
\begin{center}
\line(24,0){400}
\end{center}
\begin{abstract}
We consider dissipative periodically forced systems
and investigate cases in which having information as to how the system behaves for constant dissipation
may be used when dissipation varies in time before settling at a constant final value.
First, we consider situations where one is interested in the basins of attraction for
damping coefficients varying linearly between two given values over many different time intervals:
we outline a method to reduce the computation time required to estimate numerically the relative areas
of the basins and discuss its range of applicability.
Second, we observe that sometimes very slight changes in the time interval may produce abrupt large variations
in the relative areas of the basins of attraction of the surviving attractors: we show how comparing
the contracted phase space at a time after the final value of dissipation has been reached
with the basins of attraction corresponding to that value of constant dissipation can explain the presence of such variations.
Both procedures are illustrated by application to a pendulum with periodically oscillating support.
\end{abstract}
\begin{center}
\line(24,0){400}
\end{center}

\section{Introduction}\label{sec:1}

Consider the ordinary differential equation
\begin{equation} \label{eq:1.1}
\ddot{\theta} + f(\theta,t) + \gamma(t) \, \dot{\theta} = 0,
\end{equation}
where $\theta\in\mathbb{T}=\mathbb{R}/2\pi\mathbb{Z}$, the dots denote derivatives with respect to time $t$,
the driving force $f$ is smooth and $2\pi$-periodic with respect to its arguments
and the damping coefficient $\gamma(t) \ge 0$ depends on time.
For concreteness we shall consider explicitly the case of a pendulum with oscillating support,
which is described by an ordinary differential equation of the form
\begin{equation} \label{eq:1.2}
\ddot{\theta} + \bigl(\alpha - \beta \cos t \bigr)\sin\theta + \gamma(t) \, \dot{\theta} = 0 , \qquad
\alpha,\beta\in\mathbb{R} ,
\end{equation}
but our results apply to systems \eqref{eq:1.1} with more general forces, or even systems of the form
$\ddot x+g(x,t)+\gamma(t) \, \dot x=0$, with $x\in\mathbb{R}$, such as those considered in \cite{BDG1}.

With a few exceptions, forced systems of the form \eqref{eq:1.1} have been studied only in the case
of constant dissipation. Despite the simplicity of the model, not much is known analytically. For instance
only a finite set of attractors are expected to exist \cite{P,FGHY,RMG}, but no proof exists for this.
Furthermore, the corresponding basins of attraction are usually calculated only numerically;
see for example \cite{LH,SCVBB,AS,KG,LCFALB,BS}. To add complication,
in many physical systems the damping coefficient is not constant throughout its entire evolution.
Time-dependent dissipation has been studied in \cite{BDG1,WBG},
where the focus was mainly on the case of $\gamma(t)$ varying linearly over
some initial time span, say $t \in [0,T_0]$, after which it remains constant.
In that case, as discussed in the quoted papers, only considering the final value of the damping coefficient
does not give a correct representation of the basins of attraction
and the entire time evolution of $\gamma(t)$ must be taken into account.

Nevertheless, the analysis of a system with constant damping coefficient $\gamma$ may provide information
for the same system with damping coefficient varying in time. For fixed constant $\gamma$, consider an attractor
and denote by $A(\gamma)$ the relative area of the corresponding basin of attraction,
that is the percentage of initial data in a given sample region whose trajectories go to that attactor.
It has already been pointed out in \cite{BDG1,WBG} that, if one knows the profiles $\gamma\mapsto A(\gamma)$
describing how the relative areas
of the basins of attraction depend on $\gamma$ in the system with constant dissipation, then, when $\gamma(t)$ varies
quasi-statically (that is very slowly) from an initial value $\gamma_{i}$ towards a final value $\gamma_{f}$ over a time $T_0$,
one may be able to predict the relative areas of the basins by looking only at the profiles, and
the slower the variation of the damping coefficient, the better the prediction.
More precisely, at least in the perturbation regime, if an attractor exists for $\gamma=\gamma_{f}$,
then, when $\gamma(t)$ increases slowly towards the value $\gamma_{f}$, the relative area of the corresponding basin
of attraction is close to the value that the function $A(\gamma)$ attains for $\gamma=\gamma_{i}$.
Moreover, the larger $T_{0}$, the closer the relative area will be to this value.

In this paper, we want to discuss other cases, where knowledge of the behaviour of the system
for constant $\gamma$ may be used in the case in which dissipation changes in time.
In particular, we shall discuss the following:

\begin{enumerate}

\item Suppose that one is interested in investigating systems with $\gamma(t)$ varying linearly
from an initial value $\gamma_i$ to a final value $\gamma_{f}$, with $\gamma_i$ and $\gamma_f$ fixed,
over a time interval $[0,T_0]$, for many values of $T_0$. In Section \ref{sec:3} we outline a method
to speed up the computation of the basins of attraction.
The method utilises basins of attraction calculated for the system with constant damping coefficient
$\gamma_{f}$, and it reduces the length of time over which it is necessary for the equations
of motion to be numerically integrated. In effect, our method is to extrapolate from an observation
time $T_1 \ge T_0$ to the full evolution time $T_f$, by using pre-computed data. We will refer to this
as `the method of fast numerical computation'. This is first stated as a general method and later numerically
implemented for the pendulum with oscillating support.

\vspace{-.4cm}
\begin{figure}[H]
\centering
\subfloat[]{\includegraphics[width=0.5\textwidth]{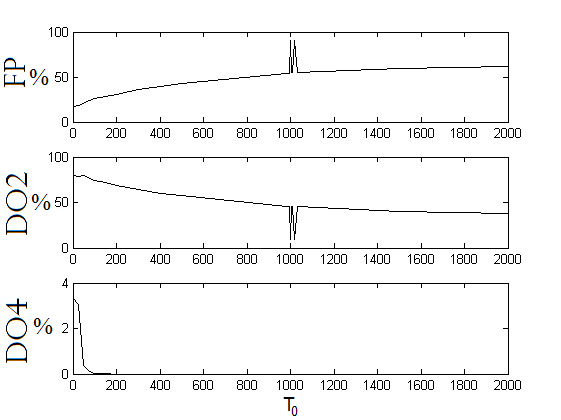}}
\subfloat[]{\includegraphics[width=0.5\textwidth]{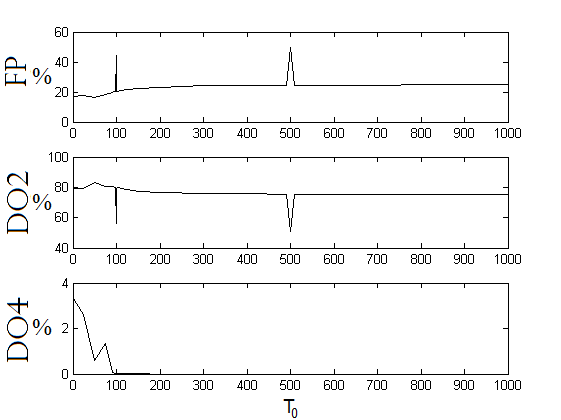}}
\caption{\small Relative areas of the basins of attraction for \eqref{eq:1.2},
with $\alpha=-0.1$, $\beta=0.545$ and $\gamma(t)$ varying linearly over time $T_0$
from (a) $\gamma_i = 0.2$ and (b) $\gamma_i = 0.23$ to $\gamma_f = 0.2725$.
FP, DO2 and DO4 denote the persisting attractors: fixed point, period-2 oscillation and
period-4 oscillation --- see Section \ref{sec:4} for details.}
\label{BasinJumps}
\end{figure}

\item In Section \ref{sec:4}, we come back to a phenomenon observed in \cite{WBG}. When the damping coefficient
increases linearly over a time $T_0$ towards a final value $\gamma_{f}$, the relative areas of the basins of attraction of the surviving
attractors most of time change smoothly as functions of $T_0$, without abrupt variations; however,
in a few cases we saw sharp jumps, concentrated in small intervals of values of $T_0$ --- see Figure \ref{BasinJumps}.
By using the same ideas as in the method of fast numerical computation, we offer an explanation of why this happens.

\end{enumerate}

The method of fast numerical computation can be described as follows.
Assume that the basins of attraction of a system of the form \eqref{eq:1.1} with constant $\gamma(t)=\gamma_f$ are known.
Of course, this information can only be obtained numerically and requires integration
of the equations of motion for a time $T_f$ large enough for the attractors to be closely
approached. Then, if one is interested in investigating
a system of the form \eqref{eq:1.1} with $\gamma(t)$ varying in time from an initial value $\gamma_i$
to the final value $\gamma_f$ over a time $T_0$ much smaller than $T_f$, the corresponding basins of attraction
can be obtained by integrating the differential equation over a time interval $[0, T_{1}]$,
with $T_{1}\ge T_0$ still much smaller than the entire time evolution $T_f$.
As a result the computation time is reduced by a factor $T_f / T_1$;
for example, with $T_1 = 48\pi$ and $T_{f}=3000$, simulations take under 12 hours
to complete, compared with nearly 10 days for the simulations using full time $T_f$.
Roughly, $T_1$ is the time needed for the trajectories to fall inside the basins of attraction at constant
$\gamma=\gamma_f$, so that their evolution from that instant onwards is known --- we refer to
Section \ref{sec:3} for more details.

This method of calculating the basins of attraction can significantly reduce computation time when
the number of values of $T_0$ to investigate is large, particularly for systems which
require long time integration or where computationally heavy integration methods are needed.
Assuming dissipation to increase with time is very natural from a physical point of view,
as an effect of wear, aging, cooling, the build up of deposits and so on. A linear increase
may seem less natural; however, in addition to being particularly well-suited
for numerical investigation, it also simulates situations where dissipation tends
to settle around some asymptotic value, which again are of physical interest; we refer
to \cite{BDG1} for further comments. Therefore, in this paper we have focused on the linear case.
Nevertheless we expect our results to apply in more general contexts, where the variation
of the damping coefficient is not linear (or not even monotonic), or the final value is reached only asymptotically.

The results produced by the method of fast numerical computation show good agreement with the results
obtained by integrating the system over the full time $T_f$ required for trajectories
to settle numerically on the persistent attractive solutions. We shall find not just that there are
only small differences between the results obtained with the method of fast numerical computation
and those of the full time integration, but also that these differences are less than the width of the 95\%
confidence interval on the results as computed by Monte Carlo simulation. However, as we shall see,
the choice of the time $T_1$ may be a delicate matter and need a few caveats.

\section{Some preliminary definitions}\label{sec:2}

In order to study the basins of attraction of a dissipative dynamical system $\dot {{\boldsymbol{x}}} =
{\boldsymbol{F}}({\boldsymbol{x}},t)$ in $\mathbb{T} \times \mathbb{R}$,
one takes initial conditions in a sample region of phase space, say $S \subset \mathbb{T} \times \mathbb{R}$,
and lets them evolve in time. As a result of dissipation, the sample region $S$ will contract with time as
trajectories move towards the persisting attractors. This leads to the following definition.

\begin{defi} \label{def:1}
Given a dissipative dynamical system and a sample region $S$, the region of phase space still
occupied by trajectories starting in $S$ at time $0$, after time $T$ has elapsed, will be called the
\emph{contracted phase space} at time $T$ and denoted by $C_{T}$.
\end{defi}

In principle $S$ is arbitrary, but in practice it is convenient to take it in such a way that (a) it contains
the (relevant) attractors and (b) one has $C_{T}\subset S$ for $T$ large enough. Thus,
if the total time over which the system evolves consists of two intervals,
a first interval $[0,T_0]$ in which $\gamma = \gamma(t)$ varies
and a second interval $[T_0,T_f]$ in which $\gamma = \gamma_f$ is constant, then whilst the first interval
takes initial conditions from the whole of $S$, the second one only receives initial conditions
from $C_{T_0}$.

It is sometimes useful to know to which points in $C_T$ points in $S$ are mapped.
This information is captured in the movement map, which is defined as follows.

\begin{defi} \label{def:2}
Fix a set $X_0$ of points ${\boldsymbol{x}}_0\in\mathbb{T}\times\mathbb{R}$ at time $t = 0$ and let $X_1$
be the set of corresponding points ${\boldsymbol{x}}_1$ that the trajectories with such initial conditions arrive at
for $t=T$. The continuous, bijective map $M_{T}$ defined by $M_{T}(X_0) = X_1$
is called a \emph{movement map} from time $0$ to time $T$.
\end{defi}

\vspace{-.5cm}
\begin{figure}[H]
\centering
\subfloat[]{\includegraphics[width=0.46\textwidth]{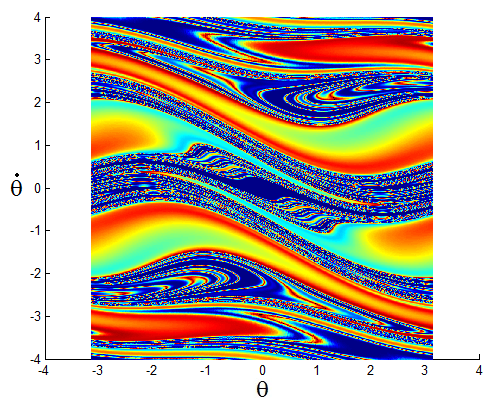}}
\subfloat[]{\includegraphics[width=0.46\textwidth]{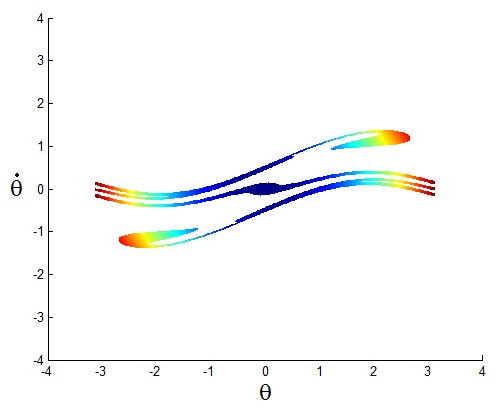}}
\caption{\small An example of movement map for \eqref{eq:1.2}, with $\alpha = -0.1$, $\beta = 0.545$
and $\gamma(t)$ growing linearly from $\gamma_i = 0.1$ to $\gamma_f = 0.2$ over time $T_0=8\pi$:
(a) sample phase space $S$ at time $t = 0$ and (b) contracted phase space  $C_{8\pi}$.
Colours are used to indicate subregions of $S$ and $C_{8\pi}$: a trajectory starting in a region of $S$
marked by a particular colour will arrive in a region of $C_{8\pi}$ with the same colour.}
\label{MovementMapExample}
\end{figure}

The contracted phase space need not (and in most cases will not) be uniformly populated. To
specify how ``dense" regions of the contracted phase space are we introduce the following definition.

\vspace{-.5cm}
\begin{figure}[H]
\centering
\subfloat[]{\includegraphics[width=0.44\textwidth]{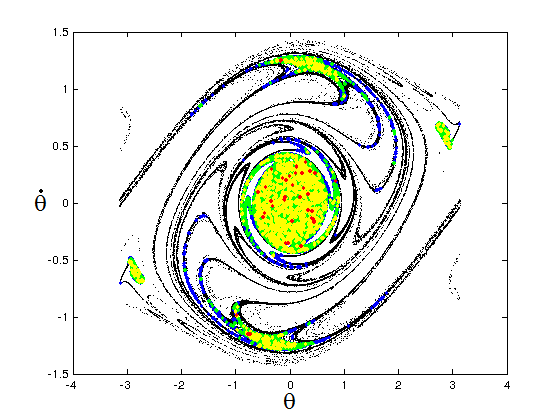}}
\subfloat[]{\includegraphics[width=0.44\textwidth]{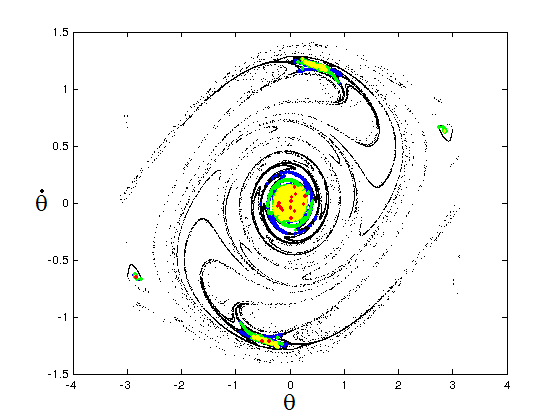}}
\caption{\small Contracted phase space for \eqref{eq:1.2},
with $\alpha = 0.5$, $\beta = 0.1$ and $\gamma(t)$ decreasing linearly from
$\gamma_i = 0.05$ to $\gamma_f = 0.02$ over a time (a) $T_0=32\pi$ and (b) $T_0=48\pi$.
The densities are coded with red being the most dense, then yellow, green, blue and black.}
\label{DensityMap}
\end{figure}

\begin{defi} \label{def:3}
Let $r>0$ be small enough. Given a dissipative dynamical system and a sample region $S$, let
$C_{T}$ be the contracted phase space at time $T$. Consider a cover of $C_{T}$
with cubes of side $r$ and denote by $X$ the set of centres of the cubes.
Fix $X_{0} \subset S$ and for ${\boldsymbol{x}}\in X$
define $\rho({\boldsymbol{x}})$ as the fraction of trajectories starting in $X_{0}$
which evolve into the cube of side $r$ and centre ${\boldsymbol{x}}$ at time $T$.
We define the \emph{density map} from time $0$ to time $T$ as the map which associates the value
$\rho({\boldsymbol{x}})$ with each ${\boldsymbol{x}} \in X$.
\end{defi}

If $X_0=S$, the function $\rho({\boldsymbol{x}})$ is expected to become $r$-independent
for $r$ small enough; for this reason we have omitted its dependence on $r$. On the other hand, if $X_0$ is a
discrete set (as in any numerical implementation), $r$ cannot be too small for the definition to make sense.

The definitions above are given explicitly in $\mathbb{T}\times\mathbb{R}$, but could be easily extended to
$\mathbb{R}^2$, or even to higher dimensional systems.
Throughout the rest of the paper we focus on systems of the form \eqref{eq:1.1}, so that
${\boldsymbol{x}}=(\theta,\dot\theta)$ and the driving force is $2\pi$-periodic in $t$.
Moreover we consider damping coefficients $\gamma(t)$ which vary linearly between two values
and become constant after a time $T_0$.
If $X_0 \subset S$ is the set of initial conditions, we take $T_{1} \ge T_{0}$ such that
$X_1 =M_{T_1}(X_0) \subset C_{T_1} \subset S$. For the system \eqref{eq:1.2},
a convenient choice for the sample region $S$ turns out to be $S=[-\pi,\pi]\times[-4,4]$.
When considering a trajectory $(\theta(t),\dot\theta(t))$,
the points $\boldsymbol{x}_0=(\theta(0),\dot{\theta}(0))$ and $\boldsymbol{x}_1=(\theta(T_1),\dot{\theta}(T_1))$
are just different points on the same trajectory and share the same path as $t \to \infty$.
Furthermore, if $T_1 = 2 N \pi$, where $N \in \mathbb{N}$, since the forcing is $2\pi$-periodic
then the solutions with initial conditions $\boldsymbol{x}_0$ and $\boldsymbol{x}_1$ at time $t=0$
move towards the same attractor if we set  $\gamma=\gamma_{f}$ from $t = 0$.
This observation will be used when choosing $T_1$ to be a multiple of $2\pi$
in the following  --- see Section \ref{sec:3.1}.

A discrete approximation to the movement map will be represented by
a $p \times m$ matrix, in which the indices $i=0,\ldots,p$ and
$j=0,\ldots, m$ correspond to the coordinates of the initial condition $(\theta_0 = -3.14 + i\Delta\theta ,\dot{\theta}_0
= -4 + j\Delta \dot{\theta})$, with $p\Delta\theta=6.28$ and $m\Delta\dot{\theta}=8$, and in each $i$,$j$-th entry
is a vector of the coordinates $(\theta_1,\dot{\theta}_1)$ for the corresponding trajectory at time $T_1$.
Colour coding the phase space we can represent the action of the movement map
as shown in Figure \ref{MovementMapExample}.
A numerical representation of a density map can be achieved by averaging trajectories onto a grid,
as done for the movement map. Different colours may be used to express how dense regions of phase space are.
In particular a density map allows us to see how the most dense regions
of the contracted phase space move with variations in $T_0$ --- see Figure \ref{DensityMap}.

\section{Fast numerical computation of the basins of attraction}\label{sec:3}

\subsection{General setting}\label{sec:3.1}

A movement map can be used to compute quickly the relative areas of basins of
attraction for systems \eqref{eq:1.1} with $\gamma(t)$ varying over a time $T_0$
before reaching a final value $\gamma_f$, with only the need to integrate
over a time $T_{1} \ge T_{0}$, provided that the basins of attraction for constant
damping coefficient $\gamma = \gamma_f$ are known. This can be achieved as follows.

First, one covers the sample region $S$ with a mesh of points,
with the requirement that the limiting solution of each point on the mesh is known
for constant $\gamma=\gamma_f$.
Of course, the basins of attraction for constant dissipation are obtained numerically, so that
which attractor any point of the mesh converges to is known only approximately;
hence, high accuracy is required at this preliminary stage.
In particular the equations of motion for initial conditions in the mesh have to be integrated
over a time interval $[0, T_f]$ sufficiently large for the attractors to be approximately reached.
Next, a set of initial conditions is chosen and for each, the equations are integrated for $t \in [0,T_1]$.
The point in phase space which the trajectory occupies at time $T_1$ is then rounded to the nearest
point on the mesh. Provided $T_1$ is chosen to be a multiple of the forcing period, we know towards
which attractor the rounded point moves asymptotically. Repeating this for all initial conditions allows one
to estimate the relative areas of the basins of attraction.

From the discussion above, it may seem natural to take $T_{1}$ to be the smallest multiple
of the forcing period $2\pi$ larger than $T_{0}$, because then a trajectory which is inside
a given basin of attraction at time $T_1$ will move towards the corresponding attractor.
However, error is inevitably introduced when approximating the coordinates of the trajectories
at time $T_1$ to the nearest point in the mesh. In systems where the basins
of attraction are sparse the slight change of coordinate during the approximation may cause
the results to be inaccurate. As a consequence it is found that in cases where $T_0$ is relatively small,
or the basins of attraction are increasingly broken and sparse, $T_1$ must be taken larger than this.
Using a finer mesh reduces the error in the approximation and should cause the results obtained
with the method of fast numerical computation described above to tend towards
those when the system is fully integrated. In practice this is
not always true, especially for systems with high sensitivity to initial conditions.
Indeed, it is possible that with a coarse mesh the approximation to the point at time $T_1$ has the same attractor,
while with a finer mesh the new, more accurate approximation moves to a different attractor.
Moreover there can be compensations which are destroyed by refining the mesh;
this is best explained by a simple example. Imagine there are two initial conditions,
one of which should tend to attractor $a_1$ and the other to attractor $a_2$,
but which, on the contrary, are predicted to tend to attractors $a_2$ and $a_1$, respectively:
in this case, the two errors cancel. Suppose we then use a finer mesh and
one of the two points is predicted correctly while the other prediction remains incorrect:
now the errors do not cancel and the predictions of the relative areas of the basins of attraction will
be in error. However, provided that this is a rare occurrence, increasing the number of points
on the mesh will tend to increase the accuracy in general.

Since basins of attraction are computed by taking a large but finite discrete set of initial conditions
(either on a mesh or uniformly and randomly distributed),
an error is always produced by using the fraction of trajectories converging towards a given attractor
to estimate the corresponding relative area.  The following statistical result gives the confidence interval
for results obtained by Monte Carlo simulation; see \cite{MMW}.
If $p$ and $\hat{p}$ are the actual and estimated probability, respectively, of landing
in a given basin of attraction, then we have
\begin{equation} \label{eq:3.1}
\hat{p} - z_{\alpha /2}\sqrt{\frac{\hat{p}(1-\hat{p})}{N}} < p < \hat{p} + z_{\alpha /2}\sqrt{\frac{\hat{p}(1-\hat{p})}{N}},
\end{equation}
where $N$ is the sample size (the number of initial conditions used for the simulation);
the variable $z_{\alpha/2}$ is the so-called $z$ value, which ensures a $(1 - \alpha)\times100\%$
confidence interval \cite{MMW}. A confidence interval of $95\%$ results in the error values shown
in Table \ref{BOAExpectedStatErr}. For example if we take $N = 10\,000$ we can see
that if an estimated basin of attraction covers 99.5\% of the  phase space,
then we can state with 95\% confidence that the actual basin of attraction will be $(99.5 \pm 0.1386)\%$.
Similarly, if the estimated basin covers 0.5\% of the phase space, the actual basin of attraction
will also be $(0.5 \pm 0.1386)\%$;  although the error is the same, the error relative to the
basin size is much greater for small basins of attraction. Thus when the basins
of attraction are small, it is necessary to use more initial conditions.

\begin{table}[H]
\centering
{\footnotesize
\rowcolors{1}{}{lightgray}
\begin{tabular}{cc|ccccccccc|}
\cline{3-11}
\multicolumn{1}{c}{} & \multicolumn{1}{c}{} & \multicolumn{9}{|c|}{\multirow{1}{*}{$N$}}\\
\cline{3-11}
\multicolumn{1}{c}{} & \multicolumn{1}{c}{} & \multicolumn{1}{|c}{10\,000} & \multicolumn{1}{c}{50\,000} & \multicolumn{1}{c}{100\,000} & \multicolumn{1}{c}{200\,000} & \multicolumn{1}{c}{300\,000} & \multicolumn{1}{c}{400\,000} & \multicolumn{1}{c}{500\,000} & \multicolumn{1}{c}{600\,000} & \multicolumn{1}{c|}{1\,000\,000} \\
\hline
\multicolumn{1}{|c|}{\multirow{12}{*}{$\hat{p} \; (\%)$}}
& \multicolumn{1}{|c|}{$0.5$} & 0.1386 & 0.0620 & 0.0438 & 0.0310 & 0.0253 & 0.0219 & 0.0196 & 0.0179 & 0.0139 \\
\multicolumn{1}{|c|}{} & \multicolumn{1}{|c|}{$1$} & 0.1950 & 0.0872 & 0.0617 & 0.0436 & 0.0356 & 0.0308 & 0.0276 & 0.0252 & 0.0195 \\
\multicolumn{1}{|c|}{} & \multicolumn{1}{|c|}{$2$} & 0.2744 & 0.1227 & 0.0868 & 0.0614 & 0.0501 & 0.0434 & 0.0388 & 0.0354 & 0.0274 \\
\multicolumn{1}{|c|}{} & \multicolumn{1}{|c|}{$3$} & 0.3344 & 0.1495 & 0.1057 & 0.0748 & 0.0610 & 0.0529 & 0.0473 & 0.0432 & 0.0334 \\
\multicolumn{1}{|c|}{} & \multicolumn{1}{|c|}{$4$} & 0.3841 & 0.1718 & 0.1215 & 0.0859 & 0.0701 & 0.0607 & 0.0543 & 0.0496 & 0.0384 \\
\multicolumn{1}{|c|}{} & \multicolumn{1}{|c|}{$5$} & 0.4272 & 0.1910 & 0.1351 & 0.0955 & 0.0780 & 0.0675 & 0.0604 & 0.0551 & 0.0427 \\
\multicolumn{1}{|c|}{} & \multicolumn{1}{|c|}{$10$} & 0.5880 & 0.2630 & 0.1859 & 0.1315 & 0.1074 & 0.0930 & 0.0832 & 0.0759 & 0.0588 \\
\multicolumn{1}{|c|}{} & \multicolumn{1}{|c|}{$15$} & 0.6999 & 0.3130 & 0.2213 & 0.1565 & 0.1278 & 0.1107 & 0.0990 & 0.0904 & 0.0700 \\
\multicolumn{1}{|c|}{} & \multicolumn{1}{|c|}{$20$} & 0.7840 & 0.3506 & 0.2479 & 0.1753 & 0.1431 & 0.1240 & 0.1109 & 0.1012 & 0.0784 \\
\multicolumn{1}{|c|}{} & \multicolumn{1}{|c|}{$30$} & 0.8982 & 0.4017 & 0.2840 & 0.2008 & 0.1640 & 0.1420 & 0.1270 & 0.1160 & 0.0898 \\
\multicolumn{1}{|c|}{} & \multicolumn{1}{|c|}{$40$} & 0.9602 & 0.4294 & 0.3036 & 0.2147 & 0.1753 & 0.1518 & 0.1358 & 0.1240 & 0.0960 \\
\multicolumn{1}{|c|}{} & \multicolumn{1}{|c|}{$50$} & 0.9800 & 0.4383 & 0.3099 & 0.2191 & 0.1789 & 0.1550 & 0.1386 & 0.1265 & 0.0980 \\ \hline
\end{tabular} }
\caption{\small The 95\% confidence interval for various relative areas, when using a Monte Carlo approach with $N$ initial
conditions, calculated from equation \eqref{eq:3.1}. This is to be interpreted as follows: if the value in the table is $x$,
then $p \in [\hat{p} - x, \hat{p} + x]$ with 95\% confidence.
An estimated relative area $\hat{p}$ has the same expected 95\%
confidence interval as $1-\hat{p}$, so it is only necessary to consider basins of attraction up to $50\%$.}
\label{BOAExpectedStatErr}
\end{table}

In the forthcoming analysis, the error in the results obtained with the method of fast numerical computation
relative to the fully integrated results will be deemed acceptable if it is less than the estimated
95\% confidence interval for the relative areas of the basins of attraction.

\subsection{Application to the pendulum with oscillating support}\label{sec:3.2}

In this section we investigate, in a concrete model, which types of basin of attraction are suitable for
our method of fast numerical computation, how to reduce the error with respect to the full integration
and how to increase the accuracy. We consider the pendulum with oscillating support, see \eqref{eq:1.2},
with $\gamma(t)$ varying linearly over a time $T_0$ from the initial
value $\gamma_i$ to the final value $\gamma_f$.
Of course, the reduction in computational time compared with
the the full time span $[0,T_f]$ is a result of the smaller integration time $T_1$.
Since $T_1 \ge T_0$, in systems where $T_0$ is large and comparable with $T_f$ this advantage is lost.
However, it was seen in \cite{WBG} that most changes to the relative areas of the basins of attraction
happen over a short initial time $T_0$, where the method is particularly effective.

The main numerical integration scheme used to test the method is {\sc matlab}'s ODE113,
which is a variable order Adams-Bashforth-Moulton scheme -- simply because {\sc matlab} offers ease in
programming compared to using a low-level language such as C. Although it is found that the integrator
ODE113 is not always reliable for this system, our aim is mainly to compare
the method of fast numerical computation with the full time integration,
rather than the accuracy of the full simulations relative to the true dynamics.
For the same reason, throughout, the relative areas of the basins of attraction are given
to 4 decimal places, despite the number of initial conditions for the full simulations producing uncertainty
in the first or second decimal place --- see Table \ref{BOAExpectedStatErr}.
As well as providing the results obtained for both fast and
full simulations with the chosen method of integration,
we also give the results obtained with
two different, more efficient integration methods, a standard Runge-Kutta integrator
and a scheme based on series expansion \cite{BDG1,WDBG}, both of which were implemented in C.
The reason is to see whether the error produced by the method of fast numerical computation is within the
difference produced by simply choosing a different numerical method and
a different selection of initial conditions or of meshes to approximate solutions at time $T_1$.

\vspace{-.3cm}
\begin{figure}[H]
\centering
\subfloat[]{\includegraphics[width=0.34\textwidth]{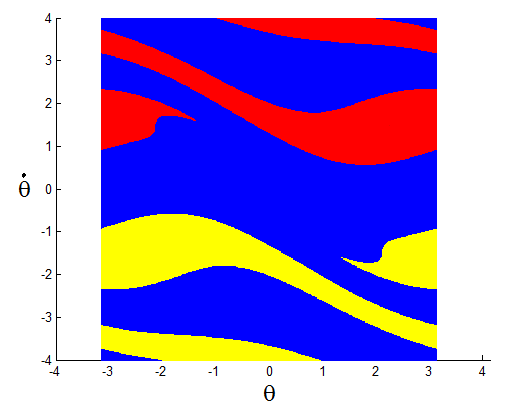}}
\subfloat[]{\includegraphics[width=0.34\textwidth]{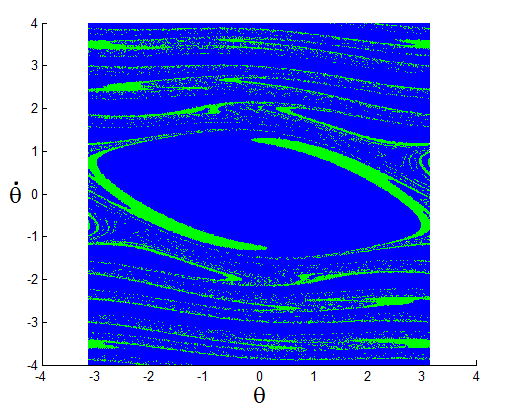}}
\subfloat[]{\includegraphics[width=0.34\textwidth]{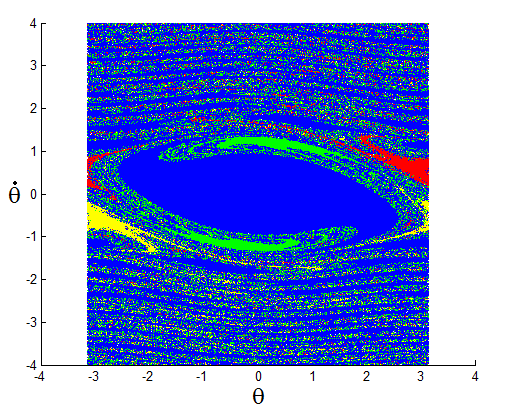}}
\caption{\small Basins of attraction for \eqref{eq:1.2} with constant
damping coefficient $\gamma(t)=\gamma_f$. In Figure (a)
$\alpha = -0.1$, $\beta = 0.545$ and $\gamma_f = 0.2$ (for clarity purposes the origin is centred
at the upwards position of the pendulum). In Figures (b) and (c)
$\alpha = 0.5$, $\beta = 0.1$ and $\gamma_f = 0.05$ and $0.02$, respectively. The
colours blue, red, yellow and green represent the basins of attraction for the attractors FP, PR,
NR and OSC, respectively.}
\label{BasinsConstantDAMM}
\end{figure}

To test the method, we consider cases where the basins of attraction for constant
values of $\gamma = \gamma_f$ become increasingly sparse. We study \eqref{eq:1.2}
with $\alpha = 0.5$, $\beta = 0.1$ and $\gamma_f = 0.05$ and $0.02$;
the corresponding basins of attraction  are shown in Figure \ref{BasinsConstantDAMM}(b)
and Figure \ref{BasinsConstantDAMM}(c), respectively.
For $\gamma=0.02$ the system exhibits four persisting attractors, namely the fixed point
$(\theta, \dot{\theta}) = (0,0)$, two rotating attractors (one positively rotating and one negatively
rotating) and an oscillatory attractor; we shall refer to them as FP, PR, NR and OSC,
respectively. For $\gamma = 0.05$, only FP and OSC persist.
Further details on these attractors can be found in \cite{BGG,WBG}.
Similar tests have also been conducted for $\alpha = -0.1$, $\beta = 0.545$
and $\gamma_f = 0.2$, with the basins of attraction as in Figure \ref{BasinsConstantDAMM}(a).
However, owing to the simple geometry of these basins of attraction, the results obtained were of little
use for studying the restrictions of the numerical application. As such the results have not been
included and hence we shall use the parameters $\alpha = 0.5$, $\beta = 0.1$ throughout the rest of this section.

\begin{figure}[H]
\centering
\subfloat[]{\includegraphics[width=0.30\textwidth]{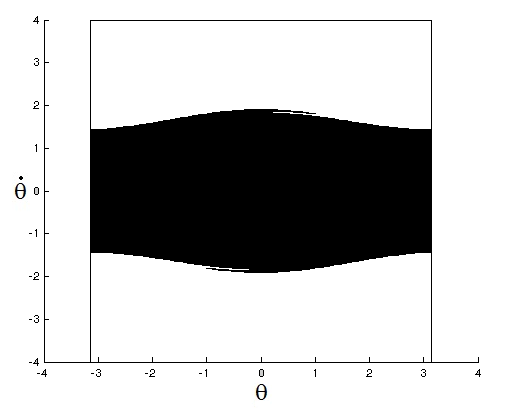}}
\hspace{1cm}
\subfloat[]{\includegraphics[width=0.30\textwidth]{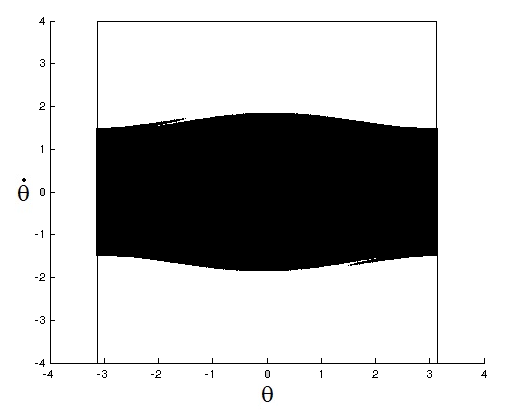}}
\vspace{-.2cm}
\caption{\small Contracted phase space for \eqref{eq:1.2} $\alpha = 0.5$, $\beta = 0.1$ and
$\gamma(t)$ varying linearly over a time $T_0= 8\pi$
(a) from $\gamma_i=0.02$ to $\gamma_f=0.05$ and (b) from $\gamma_i=0.05$ to $\gamma_f=0.02$.
The rectangle containing the shaded region marks the area $S$, from which initial conditions were taken.}
\label{PhaseSpaceRetractions8pi}
\end{figure}

In order to implement the method of fast numerical computation,
first a mesh of initial conditions must be set up and for each of these
the equations must be integrated for constant damping $\gamma_f$.
This is computationally expensive as the system is integrated over the full
time span $[0,T_f]$, required for solutions to move sufficiently close to the corresponding attractors
that they can be identified. For constant $\gamma$, this happens over times $O(1/\gamma)$
and values of $T_{f}$ between 10 and 100 times $1/\gamma$ turn out to be sufficient; for
the chosen values of $\gamma$, this yields $T_{f}$ larger than $10^{3}$.
Using the contracted phase space --- see Figure \ref{PhaseSpaceRetractions8pi} ---
for our smallest value of $T_0$ (that is $8\pi$) we can
see that it is unnecessary to cover the entire region $S$ with a mesh
of initial conditions, as roughly half of them will never be used. Only covering the region
$C_{8\pi}$ would reduce both the computation time and memory required by a factor of
roughly 2; see Figure \ref{PhaseSpaceRetractions8pi}.

\vspace{-.5cm}

\begin{figure}[H]
\centering
\subfloat[]{\includegraphics[width=0.34\textwidth]{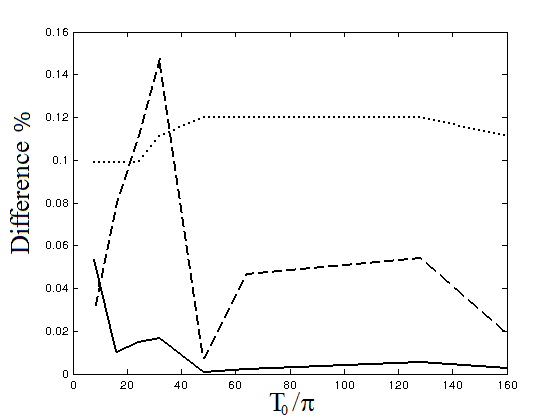}}
\hspace{1cm}
\subfloat[]{\includegraphics[width=0.34\textwidth]{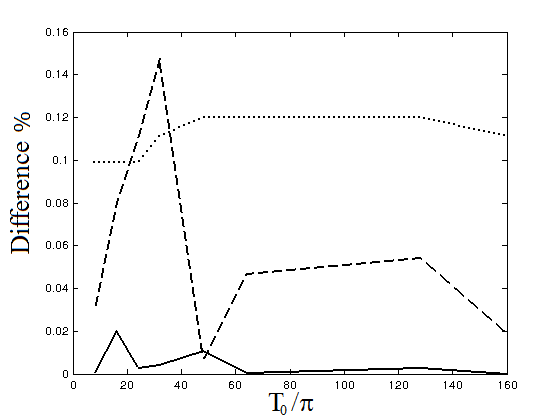}}
\caption{\small Error in the relative areas of the basins of attraction in Table \ref{DAMMError1Table}.
The full line represents the error of the method of fast numerical computation with respect to the full integration:
(a) difference between ``ODE113 Fast'' and ``ODE113",
(b) difference between ``Series Fast'' and ``Series".
The dashed line represents the difference between the estimates of ``Series'' and ``Runge-Kutta'' and the dotted line
represents the 95\% confidence interval, as calculated in Section \ref{sec:3.1}.}
\label{Difference}
\end{figure}

We first consider $\gamma(t)$ linearly evolving from $0.02$ to
$0.05$ over a time $T_0=2\pi N$, with $N\in\mathbb{N}$, and fix the integration time at $T_{1}=T_{0}$.
A mesh of roughly 500\,000 points is considered in $S$ and each trajectory at time $T_1$
is rounded to the nearest point of the mesh.
The results in Table \ref{DAMMError1Table} show that the method performs well in this instance.
Indeed, the error produced by estimating the basins of attraction with
the method of fast numerical computation rather than integrating
over the full time $T_f$ is less than the difference in results obtained by choosing a
different integrator with different initial conditions --- see Figure \ref{Difference}.
Note that the integrator ODE113 disagrees with the other integrators for $T_0 = 128\pi$ and $T_0 = 160\pi$;
however, for each integrator the results obtained with the method of fast numerical computation
are still close to the full integrations.

\begin{table}[H]
{\footnotesize
\rowcolors{1}{lightgray}{}
\begin{tabular}{c|cc!{\vrule width 0.75pt}cc!{\vrule width 0.75pt}cc!{\vrule width 0.75pt}cc!{\vrule width 0.75pt}cc|}
\cline{2-11}
&\multicolumn{10}{c|}{\multirow{1}{*}{Relative area \%}}\\
\cline{2-11}
&\multicolumn{2}{c!{\vrule width 0.75pt}}{\multirow{1}{*}{ODE113 Fast}} & \multicolumn{2}{c!{\vrule width 0.75pt}}{\multirow{1}{*}{ODE113}}
& \multicolumn{2}{c!{\vrule width 0.75pt}}{\multirow{1}{*}{Series Fast}} &
\multicolumn{2}{c!{\vrule width 0.75pt}}{\multirow{1}{*}{Series}} & \multicolumn{2}{c|}{\multirow{1}{*}{Runge-Kutta}} \\
\hline
\multicolumn{1}{|c|}{$T_0$} & \multicolumn{1}{c}{FP} & \multicolumn{1}{c|}{OSC} & \multicolumn{1}{c}{FP} & \multicolumn{1}{c|}{OSC} & \multicolumn{1}{c}{FP} & \multicolumn{1}{c|}{OSC} & \multicolumn{1}{c}{FP} & \multicolumn{1}{c|}{OSC} & \multicolumn{1}{c}{FP} & \multicolumn{1}{c|}{OSC} \\
\hline
\multicolumn{1}{|c|}{$0$} & N/A  & N/A & 85.5832 & 14.4168 & N/A  & N/A & 85.6522 & 14.3478 & 85.6744 & 14.3256 \\
\multicolumn{1}{|c|}{$8\pi$} & 86.1340 & 13.8660 & 86.0804 & 13.9196 & 86.1818 & 13.8182 & 86.1812 & 13.8188 & 86.1510 & 13.8490 \\
\multicolumn{1}{|c|}{$16\pi$} & 86.1977 & 13.8023 & 86.1874 & 13.8126 & 86.2572 & 13.7428 & 86.2370 & 13.7630 & 86.1566 & 13.8434 \\
\multicolumn{1}{|c|}{$24\pi$} & 84.0702 & 15.9298 & 84.0855 & 15.9145 & 84.1378 & 15.8622 & 84.1346 & 15.8654 & 84.0230 & 15.9770 \\
\multicolumn{1}{|c|}{$32\pi$} & 80.4622 & 19.5378 & 80.4795 & 19.5205 & 80.4374 & 19.5626 & 80.4330 & 19.5670 & 80.5804 & 19.4196 \\
\multicolumn{1}{|c|}{$48\pi$} & 75.6663 & 24.3337 & 75.6652 & 24.3348 & 75.6666 & 24.3334 & 75.6776 & 24.3224 & 75.6702 & 24.3298 \\
\multicolumn{1}{|c|}{$64\pi$} & 75.4444 & 24.5556 & 75.4472 & 24.5528 & 75.4918 & 24.5082 & 75.4912 & 24.5088 & 75.4440 & 24.5560 \\
\multicolumn{1}{|c|}{$128\pi$} & 77.3244 & 22.6756 & 77.3302 & 22.6698 & 75.8202 & 24.1798 & 75.8234 & 24.1766 & 75.8780 & 24.1220 \\
\multicolumn{1}{|c|}{$160\pi$} & 77.9221 & 22.0779 & 77.9191 & 22.0809 & 78.4282 & 21.5718 & 78.4280 & 21.5720 & 78.4092 & 21.5908 \\ \hline
\end{tabular} }
\caption{\small Relative areas of the basins of attraction for \eqref{eq:1.2} with
$\alpha = 0.5$, $\beta = 0.1$ $\gamma_i = 0.02$, $\gamma_f = 0.05$.
Initial conditions for ODE113 were taken from the same mesh (with $\Delta \theta = \Delta \dot{\theta} = 0.01$)
as used to approximate trajectories at time $T_1$, totalling 503\,829 points.
The Runge-Kutta and Series integrators each used a different set of 500\,001 random initial conditions in $S$.
The ``ODE113 Fast'' and ``Series Fast" results used the same integrator and the same initial conditions
as ``ODE113'' and ``Series", respectively.}
\label{DAMMError1Table}
\end{table}

Of course Table \ref{DAMMError1Table} allows comparison only of the relative areas of the basins of attraction,
but does not reveal to what extent the corresponding sets of points match each other. As pointed out in
Section \ref{sec:3.1}, there could be substantial compensations and, in principle, the distribution of the
basins of attraction in phase space could be wrongly described by
the method of fast numerical computation despite
the good agreement between the corresponding relative areas. However, Table \ref{PointByPointAcc}
shows that this is not the case: it is true that the percentage of points assigned by
the method of fast numerical computation to
the wrong basin is larger than suggested by Table \ref{DAMMError1Table},
but still is not large, and also tends to decrease with increasing $T_0$ and hence $T_1$.

\begin{figure}[H]
\centering
\subfloat[]{\includegraphics[width=0.34\textwidth]{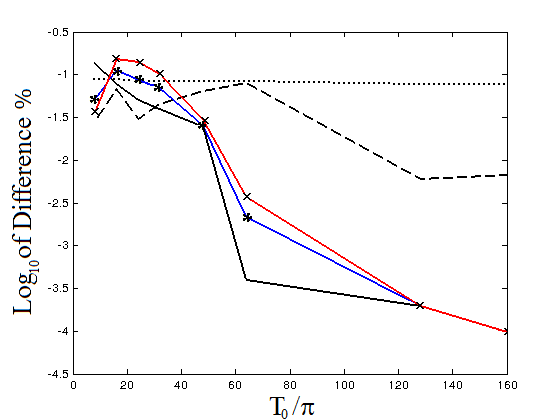}}
\hspace{1cm}
\subfloat[]{\includegraphics[width=0.34\textwidth]{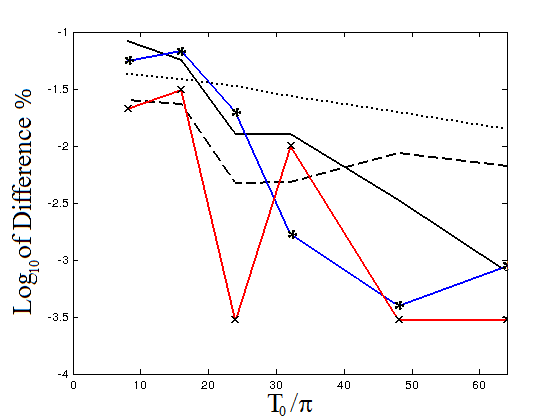}}\\
\subfloat[]{\includegraphics[width=0.34\textwidth]{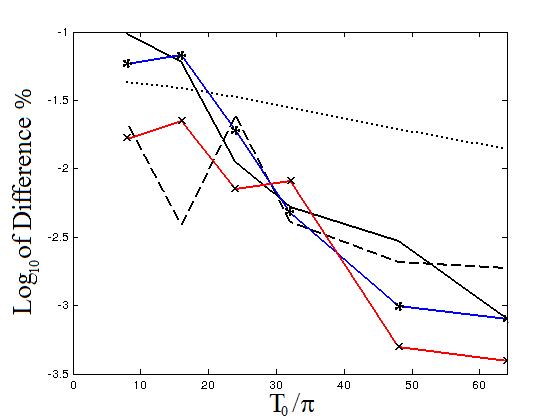}}
\hspace{1cm}
\subfloat[]{\includegraphics[width=0.34\textwidth]{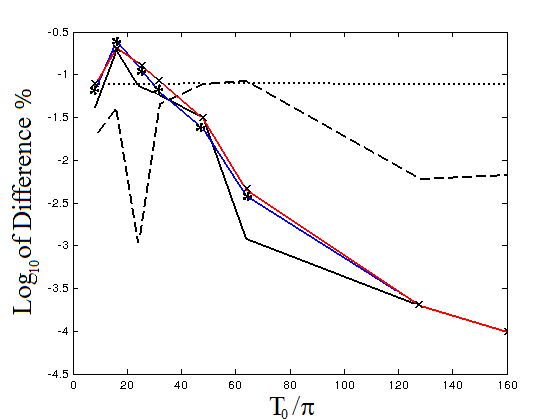}}
\caption{\small Error in the relative areas of the basins of attraction in Table \ref{DAMMError3Table}:
Figures (a) through to (d) show the differences for FP, PR, NR and OSC, respectively, with
(b) and (c) only going as far as $T_0=64\pi$ since the rotating solutions disappear after this.
The plain, starred and crossed full lines show the difference between ``Fast 1", ``Fast 2"
and ``Fast 3", respectively, and ``ODE113''. The dashed line shows the difference
between the estimates using ``Series'' and ``Runge-Kutta''. The dotted line shows
the 95\% confidence interval.}
\label{Difference-1000000}
\end{figure}

This may still be regarded as a very simple case, as only
two attractors coexist and the boundaries of their basins of attraction
for constant $\gamma_f$ still have relatively simple geometry.
We now consider \eqref{eq:1.2},
with the same values for $\alpha$ and $\beta$ as in the previous case,
but with $\gamma(t)$ linearly decreasing from $\gamma_i=0.05$ to $\gamma_{f}=0.02$.
We first ran simulations with roughly 500\,000 initial conditions and two meshes with increments
$(\Delta \theta,\Delta \dot{\theta})=(0.01,0.01)$ and $(0.005,0.005)$, respectively;
at time $T_{1}=T_0$ each trajectory was rounded to the nearest point of the mesh.
Afterwards we considered roughly 1\,000\,000 initial conditions and three meshes with increments
$(\Delta \theta,\Delta \dot{\theta})=(0.01,0.01)$, $(0.01,0.005)$ and $(0.005,0.005)$, respectively;
once more we chose $T_1=T_0$.
The corresponding results are reported in Table \ref{DAMMError3Table}.
Figure \ref{Difference-1000000} illustrates that, as previously,
the difference between the fast numerical computations and the full time integrations
is comparable to --- if not smaller than --- both the difference between the results obtained with
the other two integrators and the 95\% confidence interval.

\vspace{-.3cm}

\vspace{-.5cm}
\begin{table}[H]
\centering
{\footnotesize
\rowcolors{1}{lightgray}{}
\begin{tabular}{c|cc|}
\cline{2-3}
&\multicolumn{2}{c|}{\multirow{1}{*}{Accuracy of the method \%}}\\
\hline
\multicolumn{1}{|c|}{$T_0$} & \multicolumn{1}{c}{ODE113 Fast} & \multicolumn{1}{c|}{Series Fast}\\
\hline
\multicolumn{1}{|c|}{$8\pi$} & 96.6004 & 97.1218 \\
\multicolumn{1}{|c|}{$16\pi$} & 99.0513 & 99.1206 \\
\multicolumn{1}{|c|}{$24\pi$} & 99.6497 & 99.6984 \\
\multicolumn{1}{|c|}{$32\pi$} & 99.8577 & 99.8860 \\
\multicolumn{1}{|c|}{$48\pi$} & 99.9464 & 99.9686 \\
\multicolumn{1}{|c|}{$64\pi$} & 99.9678 & 99.9922 \\
\multicolumn{1}{|c|}{$128\pi$} & 99.7600 & 99.9999 \\
\multicolumn{1}{|c|}{$160\pi$} & 99.5254 & 99.99998 \\ \hline
\end{tabular} }
\vspace{-.1cm}
\caption{\small The percentages of initial conditions for which the limiting behaviour is
correctly described for the simulations in Table \ref{DAMMError1Table} using
the method of fast numerical computation instead of the full integration.}
\label{PointByPointAcc}
\end{table}

As expected, the difference between the results obtained using different integrators
reduces noticeably when increasing the number of initial conditions. Contrary to this,
neither increasing the number of initial conditions nor using a finer mesh
in general significantly improves the error relative to the full time integration ---
see Table \ref{PointByPointAcc2}. Slight improvements are only obtained for the attractors
with smaller basin of attraction (that is the two rotating attractors), where the error
becomes comparable to that created by using different numerical integrators and sets of initial conditions.

\begin{figure}[H]
\centering
\subfloat[]{\includegraphics[width=0.34\textwidth]{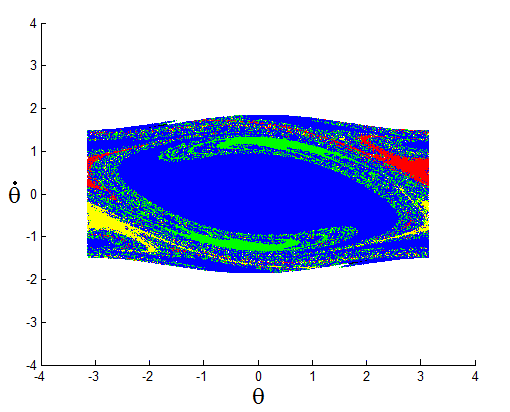}}
\subfloat[]{\includegraphics[width=0.34\textwidth]{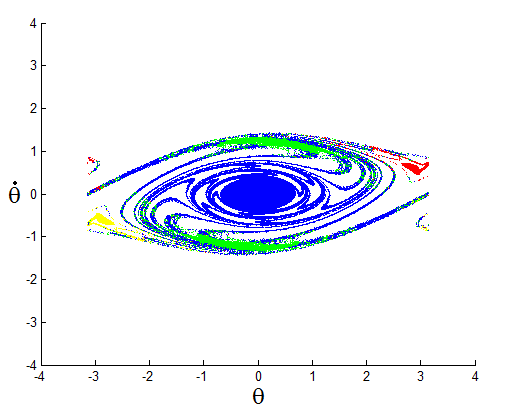}}
\subfloat[]{\includegraphics[width=0.34\textwidth]{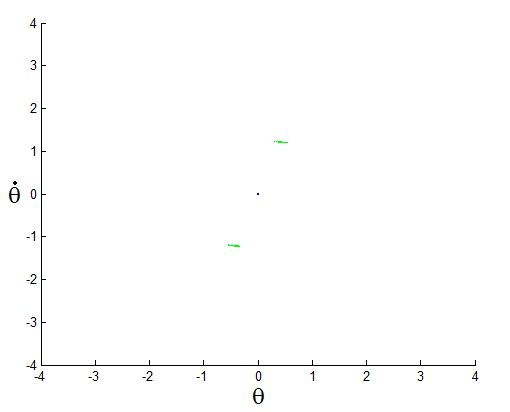}}
\caption{\small Images of the contracted phase space $C_{T_1}$ superimposed on top of the basins
of attraction for constant $\gamma_f = 0.02$. The values of $T_1$ are (a) $8\pi$, (b) $32\pi$ and
(c) $160\pi$.}
\label{From005to002RetractionFig}
\end{figure}

It is also apparent that the larger the value of $T_0$ (and thus also $T_1$), the more reliable the method.
Table \ref{PointByPointAcc2} shows that for $T_0 \ge 48\pi$
the results are correct for more than 99\% of initial conditions.
The reason behind this can be seen by plotting the contracted phase space $C_{T_1}$ for different
values of $T_1$ and superimposing it on the basins of attraction for constant
$\gamma = \gamma_f= 0.02$; see Figure \ref{From005to002RetractionFig}.
When $T_1$ is larger, $C_{T_1}$
occupies regions of the basin of attraction which are
deep inside the cores surrounding the attractors; thus, they are distant from the boundaries and
hence less sensitive to slight variations in initial conditions. In turn, the error created by the process
of approximation onto the mesh is less significant. For the same reason, as
the results in Table \ref{DAMMError3Table} corresponding to $T_0 = 128\pi$ and $160\pi$
show, decreasing the spacing of points in the mesh does not improve the results for $T_0$ large.

\vspace{.2cm}
\begin{table}[H]
\centering
\footnotesize{
\rowcolors{1}{lightgray}{}
\begin{tabular}{c|ccccccc|}
\cline{2-8}
&&\multicolumn{6}{c|}{\multirow{1}{*}{Relative area \%}}\\
\hline
\multicolumn{1}{|c|}{$T_0$} & \multicolumn{1}{c}{Attractor} & \multicolumn{1}{c}{Fast 1} & \multicolumn{1}{c}{Fast 2} & \multicolumn{1}{c}{Fast 3} & \multicolumn{1}{c}{ODE113} & \multicolumn{1}{c}{Series}
& \multicolumn{1}{c|}{Runge-Kutta} \\
\hline
\multicolumn{1}{|c|}{\multirow{4}{*}{$0$}} & \multicolumn{1}{c|}{FP} & N/A & N/A & N/A & 71.8847 & 71.2994 & 71.3856 \\
\multicolumn{1}{|c|}{}& \multicolumn{1}{c|}{PR} & N/A & N/A & N/A & 4.6098 & 4.8996 & 4.8844 \\
\multicolumn{1}{|c|}{}& \multicolumn{1}{c|}{NR} & N/A & N/A & N/A & 4.6149 & 4.8800 & 4.8916 \\
\multicolumn{1}{|c|}{}& \multicolumn{1}{c|}{OSC} & N/A & N/A & N/A & 18.8906 & 18.9210 & 18.8384 \\ \hline
\multicolumn{1}{|c|}{\multirow{4}{*}{$8\pi$}} & \multicolumn{1}{c|}{FP} & 71.9491 & 72.0382 & 72.1258 & 72.0892 & 71.4668 & 71.4953 \\
\multicolumn{1}{|c|}{}& \multicolumn{1}{c|}{PR} & 4.6425 & 4.6150 & 4.5799 & 4.5587 & 4.8784 & 4.8529 \\
\multicolumn{1}{|c|}{}& \multicolumn{1}{c|}{NR} & 4.6513 & 4.6132 & 4.5706 & 4.5540 & 4.8703 & 4.8484 \\
\multicolumn{1}{|c|}{}& \multicolumn{1}{c|}{OSC} & 18.7572 & 18.7336 & 18.7237 & 18.7982 & 18.7845 & 18.8034 \\ \hline
\multicolumn{1}{|c|}{\multirow{4}{*}{$16\pi$}} & \multicolumn{1}{c|}{FP} & 73.2636 & 73.3001 & 73.3398 & 73.1836 & 72.6196 & 72.6889 \\
\multicolumn{1}{|c|}{}& \multicolumn{1}{c|}{PR} & 4.2390 & 4.2512 & 4.2130 & 4.1818 & 4.4743 & 4.4508 \\
\multicolumn{1}{|c|}{}& \multicolumn{1}{c|}{NR} & 4.2441 & 4.2515 & 4.2054 & 4.1828 & 4.4691 & 4.4652 \\
\multicolumn{1}{|c|}{}& \multicolumn{1}{c|}{OSC} & 18.2533 & 18.1972 & 18.2418 & 18.4518 & 18.4370 & 18.3951 \\ \hline
\multicolumn{1}{|c|}{\multirow{4}{*}{$24\pi$}} & \multicolumn{1}{c|}{FP} & 76.1441 & 76.1815 & 76.2378 & 76.0930 & 76.0347 & 76.0037 \\
\multicolumn{1}{|c|}{}& \multicolumn{1}{c|}{PR} & 3.1458 & 3.1531 & 3.1333 & 3.1330 & 3.1620 & 3.1668 \\
\multicolumn{1}{|c|}{}& \multicolumn{1}{c|}{NR} & 3.1466 & 3.1548 & 3.1281 & 3.1353 & 3.1826 & 3.2077 \\
\multicolumn{1}{|c|}{}& \multicolumn{1}{c|}{OSC} & 17.5635 & 17.5105 & 17.5008 & 17.6387 & 17.6207 & 17.6218 \\ \hline
\multicolumn{1}{|c|}{\multirow{4}{*}{$32\pi$}} & \multicolumn{1}{c|}{FP} & 76.9920 & 77.0224 & 77.0556 & 76.9517 & 76.9123 & 76.8659 \\
\multicolumn{1}{|c|}{}& \multicolumn{1}{c|}{PR} & 2.1505 & 2.1358 & 2.1274 & 2.1375 & 2.0973 & 2.1022 \\
\multicolumn{1}{|c|}{}& \multicolumn{1}{c|}{NR} & 2.1416 & 2.1315 & 2.1281 & 2.1363 & 2.1171 & 2.1130 \\
\multicolumn{1}{|c|}{}& \multicolumn{1}{c|}{OSC} & 18.7158 & 18.7103 & 18.6888 & 18.7745 & 18.8733 & 18.9189 \\ \hline
\multicolumn{1}{|c|}{\multirow{4}{*}{$48\pi$}} & \multicolumn{1}{c|}{FP} & 76.3075 & 76.3077 & 76.3121 & 76.2824 & 76.2565 & 76.3226 \\
\multicolumn{1}{|c|}{}& \multicolumn{1}{c|}{PR} & 0.9210 & 0.9172 & 0.9179 & 0.9176 & 0.8702 & 0.8789 \\
\multicolumn{1}{|c|}{}& \multicolumn{1}{c|}{NR} & 0.9209 & 0.9169 & 0.9184 & 0.9179 & 0.8750 & 0.8771 \\
\multicolumn{1}{|c|}{}& \multicolumn{1}{c|}{OSC} & 21.8505 & 21.8582 & 21.8515 & 21.8820 & 21.9983 & 21.9189 \\ \hline
\multicolumn{1}{|c|}{\multirow{4}{*}{$64\pi$}} & \multicolumn{1}{c|}{FP} & 77.6648 & 77.6674 & 77.6690 & 77.6652 & 77.6993 & 77.7806 \\
\multicolumn{1}{|c|}{}& \multicolumn{1}{c|}{PR} & 0.2970 & 0.2971 & 0.2965 & 0.2962 & 0.2752 & 0.2819 \\
\multicolumn{1}{|c|}{}& \multicolumn{1}{c|}{NR} & 0.2972 & 0.2972 & 0.2968 & 0.2964 & 0.2758 & 0.2739 \\
\multicolumn{1}{|c|}{}& \multicolumn{1}{c|}{OSC} & 21.7410 & 21.7383 & 21.7377 & 21.7422 & 21.7497 & 21.6636 \\ \hline
\multicolumn{1}{|c|}{\multirow{2}{*}{$128\pi$}} & \multicolumn{1}{c|}{FP} & 81.0690 & 81.0690 & 81.0690 & 81.0688 & 81.1812 & 81.1750 \\
\multicolumn{1}{|c|}{}& \multicolumn{1}{c|}{OSC} & 18.9310 & 18.9310 & 18.9310 & 18.9312 & 18.8188 & 18.8250 \\ \hline
\multicolumn{1}{|c|}{\multirow{2}{*}{$160\pi$}} & \multicolumn{1}{c|}{FP} & 81.8492 & 81.8492 & 81.8492 & 81.8491 & 81.9317 & 81.9386 \\
\multicolumn{1}{|c|}{}& \multicolumn{1}{c|}{OSC} & 18.1508 & 18.1508 & 18.1508 & 18.1509 & 18.0683 & 18.0614 \\ \hline
\end{tabular} }
\vspace{-0.1cm}
\caption{\small Results for $\gamma(t)$ varying from $\gamma_i = 0.05$ to $\gamma_f = 0.02$ over
times $T_0$. The Runge-Kutta and Series integrators used different sets of 1\,000\,001 random
initial conditions in $S$, while ODE113 used a mesh of initial conditions in $S$
with increments $(\Delta \theta, \Delta \dot{\theta}) = (0.01,0.005)$.
The meshes in ``Fast 1'', ``Fast 2'' and ``Fast 3''
have increments $(\Delta \theta, \Delta \dot{\theta}) = (0.01,0.01)$, $(0.01,0.005)$ and $(0.005,0.005)$, respectively.}
\label{DAMMError3Table}
\end{table}

\vspace{-0.5cm}
\begin{table}[H]
\centering
{\footnotesize
\rowcolors{1}{lightgray}{}
\begin{tabular}{c|cc!{\vrule width 0.75pt}ccc|}
\cline{2-6}
&\multicolumn{5}{c|}{\multirow{1}{*}{Accuracy of the method \%}}\\
\cline{2-6}
&\multicolumn{2}{c!{\vrule width 0.75pt}}{\multirow{1}{*}{$N \approx \text{500\,000}$}} & \multicolumn{3}{c!{\vrule width 0.75pt}}{\multirow{1}{*}{$N \approx \text{1\,000\,000}$}} \\
\hline
\multicolumn{1}{|c|}{$T_0$}& \multicolumn{1}{c}{Mesh 1} & \multicolumn{1}{c|}{Mesh 3} & \multicolumn{1}{c}{Mesh 1} & \multicolumn{1}{c}{Mesh 2} & \multicolumn{1}{c|}{Mesh 3} \\
\hline
\multicolumn{1}{|c|}{$8\pi$} & 77.1284 & 77.1986 & 76.5894 & 75.8780 & 78.0587 \\
\multicolumn{1}{|c|}{$16\pi$} & 87.9880 & 87.9715 & 88.0697 & 86.8358 & 88.7351 \\
\multicolumn{1}{|c|}{$24\pi$} & 95.0983 & 95.4209 & 95.0905 & 93.6584 & 95.4112 \\
\multicolumn{1}{|c|}{$32\pi$} & 97.4872 & 97.7471 & 97.5068 & 96.3252 & 97.7284 \\
\multicolumn{1}{|c|}{$48\pi$} & 99.6298 & 99.6542 & 99.6308 & 99.3272 & 99.6504 \\
\multicolumn{1}{|c|}{$64\pi$} & 99.9637 & 99.9637 & 99.9585 & 99.9262 & 99.9586 \\
\multicolumn{1}{|c|}{$128\pi$} & 99.9996 & 99.9996 & 99.9996 & 99.9996 & 99.9996 \\
\multicolumn{1}{|c|}{$160\pi$} & 99.9998 & 99.9998 & 99.9999 & 99.9999 & 99.9999 \\ \hline
\end{tabular} }
\vspace{-0.1cm}
\caption{\small The percentages of initial conditions for which the limiting behaviour is
correctly described for the simulations in Table \ref{DAMMError3Table} and the analogous ones
with 500\,000 initial conditions. Mesh 1, Mesh 2 and Mesh 3 correspond to increments
$(\Delta \theta, \Delta \dot{\theta}) = (0.01,0.01)$, $(0.01,0.005)$ and $(0.005,0.005)$, respectively.}
\label{PointByPointAcc2}
\end{table}

\vspace{-0.4cm}
\begin{table}[H]
\centering
{\footnotesize
\rowcolors{1}{lightgray}{}
\begin{tabular}{c|cccc!{\vrule width 0.75pt}cccc|}
\cline{2-9}
&\multicolumn{8}{c|}{\multirow{1}{*}{Relative area \%}}\\
\cline{2-9}
&\multicolumn{4}{c!{\vrule width 0.75pt}}{\multirow{1}{*}{Fast 1}} & \multicolumn{4}{c!{\vrule width 0.75pt}}{\multirow{1}{*}{Fast 3}} \\
\hline
\multicolumn{1}{|c|}{$T_0$} & \multicolumn{1}{c}{FP} & \multicolumn{1}{c}{PR} & \multicolumn{1}{c}{NR} & \multicolumn{1}{c|}{OSC} & \multicolumn{1}{c}{FP} & \multicolumn{1}{c}{PR} & \multicolumn{1}{c}{NR} & \multicolumn{1}{c|}{OSC} \\
\hline
\multicolumn{1}{|c|}{$8\pi$} & 72.1750 & 4.5484 & 4.5771 & 18.6995 & 72.1704 & 4.5795 & 4.5661 & 18.6838 \\
\multicolumn{1}{|c|}{$16\pi$} & 73.2071 & 4.1834 & 4.2040 & 18.4054 & 73.2283 & 4.1922 & 4.1835 & 18.3960 \\
\multicolumn{1}{|c|}{$24\pi$} & 76.1309 & 3.1296 & 3.1390 & 17.6005 & 76.1392 & 3.1359 & 3.1327 & 17.5922  \\
\multicolumn{1}{|c|}{$32\pi$} & 76.9679 & 2.1347 & 2.1343 & 18.7631 & 76.9765 & 2.1371 & 2.1350 & 18.7514 \\ \hline
\end{tabular} }
\vspace{-0.1cm}
\caption{\small Results for $\gamma(t)$ varying from $\gamma_i = 0.05$ to $\gamma_f = 0.02$ over times $T_0$,
with integration time $T_1 = 48\pi$. The meshes for ``Fast 1'' and ``Fast 3''
have increments of $(0.01,0.01)$ and $(0.005,0.005)$, respectively.}
\label{ImprovedErrorsTable}
\end{table}

\null\vspace{-.8cm}

\vspace{-.4cm}

\begin{table}[H]
\centering
{\footnotesize
\rowcolors{1}{lightgray}{}
\begin{tabular}{c|c|c|}
\cline{2-3}
&\multicolumn{2}{c|}{\multirow{1}{*}{Accuracy of the method \%}}\\
\cline{2-3}
&\multicolumn{1}{c!{\vrule width 0.75pt}}{\multirow{1}{*}{Fast 1}} & \multicolumn{1}{c!{\vrule width 0.75pt}}{\multirow{1}{*}{Fast 3}} \\
\hline
\multicolumn{1}{|c|}{$8\pi$} & 94.8805 & 95.1295 \\
\multicolumn{1}{|c|}{$16\pi$} & 96.2629 & 96.4448 \\
\multicolumn{1}{|c|}{$24\pi$} & 98.0571 & 98.1563 \\
\multicolumn{1}{|c|}{$32\pi$} & 98.7308 & 98.8242 \\ \hline
\end{tabular} }
\caption{\small The percentages of initial conditions for which the limiting behaviour is correctly
described for the simulations in Table \ref{ImprovedErrorsTable} using
an integration time $T_1=48\pi$ instead of $T_1=T_0$.
The results are less than 99\% accurate because $\gamma_i > \gamma_f$ and hence $S$ contracts
more slowly at $\gamma = \gamma_f$, i.e. in the region $[T_0, T_1]$.}
\label{ImprovedDAMMaccuracy}
\end{table}

\vspace{-.6cm}
\begin{figure}[H]
\centering
\subfloat[]{\includegraphics[width=0.34\textwidth]{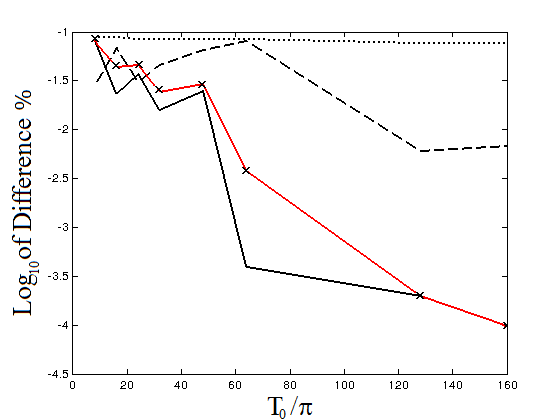}}
\hspace{1cm}
\subfloat[]{\includegraphics[width=0.34\textwidth]{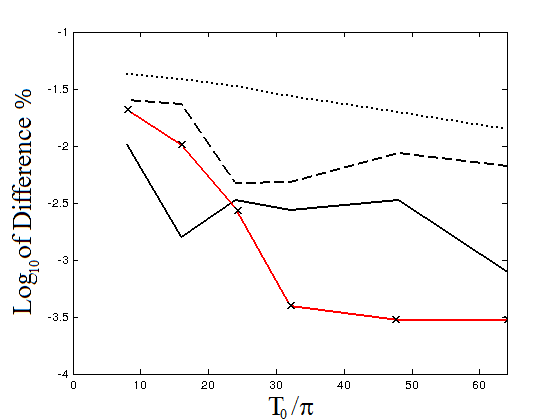}}\\
\subfloat[]{\includegraphics[width=0.34\textwidth]{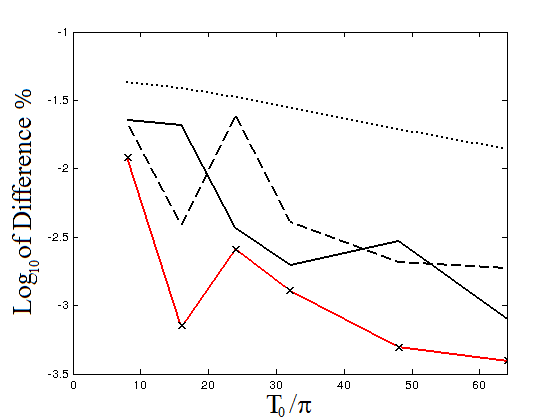}}
\hspace{1cm}
\subfloat[]{\includegraphics[width=0.34\textwidth]{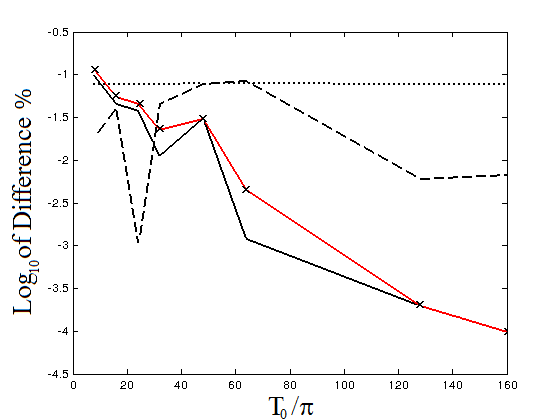}}
\caption{\small
The same as Figure \ref{Difference-1000000} with an integration time $T_{1}=\max\{48\pi,T_0\}$.
The plain and crossed full lines show the differences between
``Fast 1" and ``Fast 3", respectively, and ``ODE113".}
\label{ImprovedErrorsFigure}
\end{figure}

All this suggests that setting a minimum value for the integration time $T_1$
increases the accuracy of the method of fast numerical computation.
Doing so also results in an increase in the computational time if $T_{1} \gg T_{0}$.
However it is still quicker than integrating over the full time $T_{f}$.
Moreover the larger $T_1$, the less of phase space the mesh is required to cover,
so reducing the points in the mesh.
Tables \ref{ImprovedErrorsTable} and \ref{ImprovedDAMMaccuracy} give
the results obtained by integrating solutions over a time $T_1 =48\pi$ instead of $T_1=T_0$,
to be compared with those in Table \ref{DAMMError3Table} and \ref{PointByPointAcc2}, respectively.
Figure \ref{ImprovedErrorsFigure} shows that the error has dramatically decreased,
falling below the differences between the other two integrators, even when a coarse mesh is used
(to reduce further the overall integration time).

\section{Jumps in the relative areas of the basins of attraction} \label{sec:4}

\begin{figure}[H]
\centering
\subfloat[]{\includegraphics[width=0.34\textwidth]{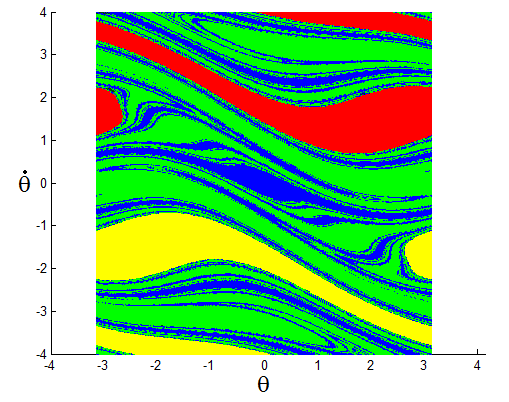}}
\subfloat[]{\includegraphics[width=0.34\textwidth]{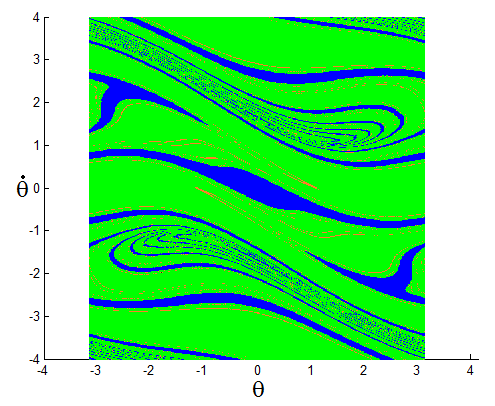}}
\vspace{-.1cm}
\caption{\small Basins of attraction for \eqref{eq:1.2} with $\alpha=-0.1$, $\beta=0.545$ and
(a) $\gamma = 0.23$ and (b) $0.2725$. The basins of attraction
are as follows: Blue = FP, Green = DO2, Red=PR, Yellow=NR and Orange = DO4.}
\label{BasinsConst}
\end{figure}

\begin{figure}[H]
\centering
\subfloat[]{\includegraphics[width=0.34\textwidth]{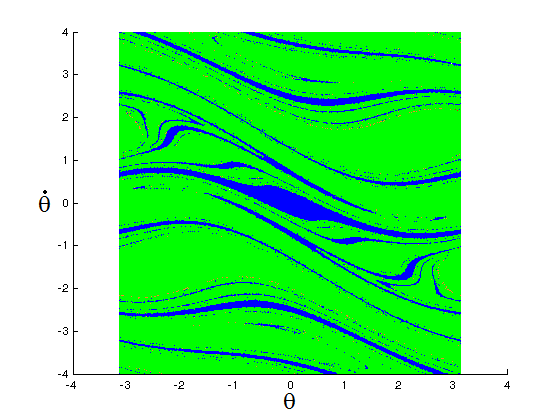}}
\subfloat[]{\includegraphics[width=0.34\textwidth]{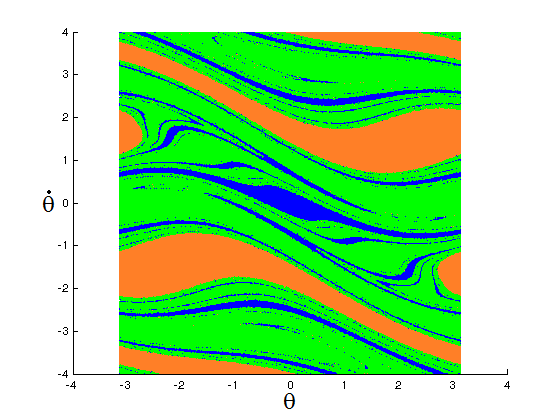}}
\subfloat[]{\includegraphics[width=0.34\textwidth]{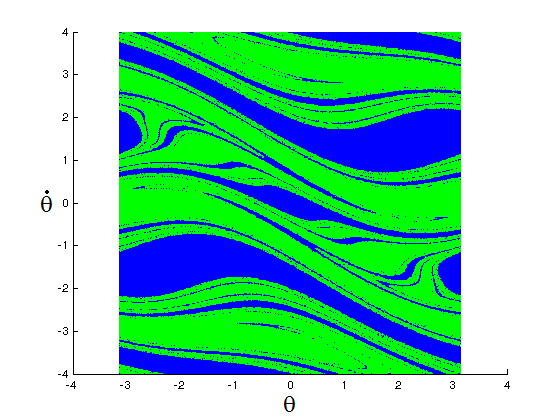}}
\vspace{-.1cm}
\caption{\small Basins of attraction for \eqref{eq:1.2} with $\alpha=-0.1$, $\beta=0.545$ and
$\gamma(t)$ increasing linearly from $\gamma_i = 0.23$ to
$\gamma_f=0.2725$ over a time (a) $T_0=93.6$, (b) $94.0$ and (c) $99.2$. The colours
are as in Figure \ref{BasinsConst}.}
\label{BasinsNotConst}
\end{figure}

We consider \eqref{eq:1.2} with $\alpha = -0.1$, $\beta = 0.545$ and $\gamma(t)$ varying from $\gamma_i = 0.23$
to $\gamma_f = 0.2725$. For $\gamma=\gamma_{i}$ the system admits four attractors: the fixed point (FP),
the period-1 positively and negatively rotating solutions (PR/NR) and the period-2 oscillations (DO2) ---
we refer to \cite{WBG} for more details. For
$\gamma=\gamma_{f}$ only FP and DO2 survive and a new attractor appears: the period-4 oscillation
(DO4). The corresponding basins of attraction are shown in Figure \ref{BasinsConst};
in Figure \ref{BasinsNotConst} the basins of attraction for $\gamma(t)$ varying in time from $\gamma_i=0.23$
to $\gamma_f=0.2725$ are represented for some values of $T_0$.

\begin{figure}[H]
\centering
\subfloat[]{\includegraphics[width=0.52\textwidth]{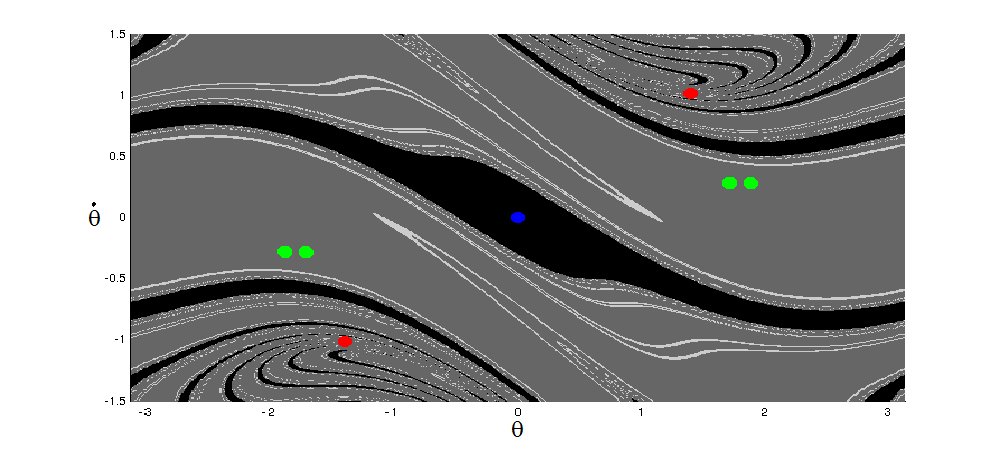}}
\subfloat[]{\includegraphics[width=0.52\textwidth]{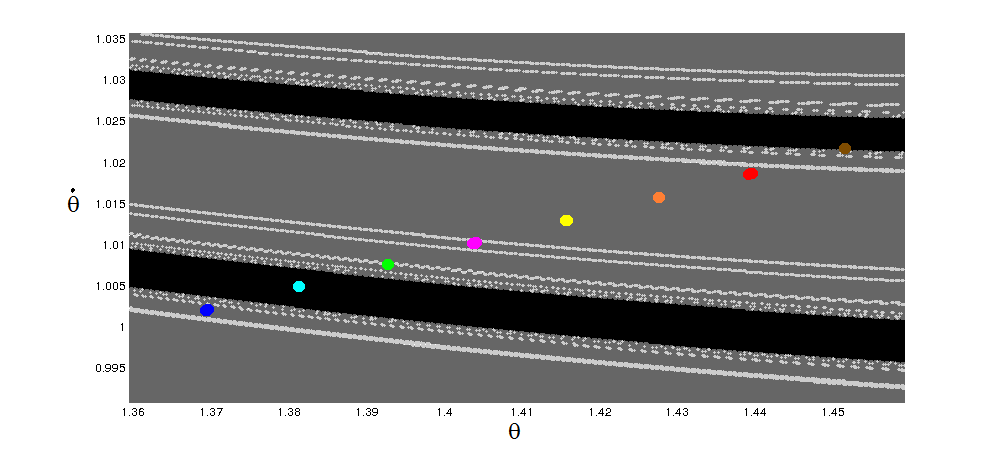}}
\caption{\small Clusters in phase space superimposed onto basins of attraction
for constant $\gamma=0.2725$:
(a) whole sample region  with clusters FP, DO2 and U in Blue, Green and Red, respectively,
and (b) magnification of the region
containing the clusters U in the upper half-plane (here the colour code is:
Blue=92, Cyan=93, Green=94, Magenta=95, Yellow=96, Orange=97, Red=98, Brown=99).
The basins are colour coded with Black for FP, Dark Grey for DO2 and Light Grey for DO4.
The snapshots in time are taken at time $T_{1}=32\pi$. The coloured
dots have been drawn larger than the actual sizes of the clusters for clarity purposes.}
\label{DenseRegionsAll}
\end{figure}

In \cite{WBG}, for $\gamma(t)$ varying from $\gamma_i = 0.23$ to $\gamma_f = 0.2725$ over a time $T_{0}$,
we noticed large transitions in the basins of attraction for $T_0 = 99$, 100, 500; see Figure \ref{BasinJumps}(b).
To understand why a jump appears, say, for $T_0=99$, we can reason as follows.
We increase $T_{0}$ in steps of one from 92 to 99. For all such values the appropriate density map
shows that, after a time $T_1 = 32\pi > T_0$, over 90\% of the trajectories are clustered in very small regions:
this is an effect of the dissipation being rather large.
Indeed, all trajectories start tending to the attractors existing at $\gamma = 0.23$.
However, with the dissipation changing, also the attractors change quasi-statically.
The attractors FP and DO2 still persist at $\gamma=0.2725$. Most of the trajectories moving towards them
form clusters at time $T_1$ close to the points at which the attractors cross the plane at $t=T_1$;
we denote by FP and DO2 such clusters, according to the attractor they are approaching.
Instead, the attractors PR and NR disappear when the damping coefficient reaches
the value $\gamma \approx 0.269$ \cite{WBG}. The trajectories that were moving towards them
form clusters that we denote by U. Thus, at time $T_1$, there are seven clusters of points:
two clusters $U$, one cluster FP and four clusters DO2 (since up to $\gamma \approx 0.27$ there are
two period-2 solutions) --- see Figure \ref{DenseRegionsAll}.
The cluster FP and each cluster DO2 correspond roughly to 20\% and
11-13\% of initial conditions, respectively, while the two clusters U each represent roughly 12\% of the initial conditions
(approximately the relative areas of the basins of attraction for the rotating solutions at $\gamma = 0.23$).

\vspace{-.5cm}

\vspace{-1.0cm}

Since both attractors FP and DO2 remain stable for all values
of $\gamma\in[\gamma_i , \gamma_f]$, it is unlikely that
the clusters  labelled as FP and DO2  will change their limiting solutions. This is confirmed by the fact that the clusters
FP and DO2 are well inside the basins of attraction of the corresponding attractors FP and DO2,
respectively, for all values of $T_0$.
On the contrary the clusters U move slightly from right to left in the upper half-plane
and from left to right in the lower half-plane as $T_0$ increases.
In doing so, they cross the boundaries of the basins of attraction for constant $\gamma = 0.2725$.

\begin{figure}[H]
\centering
\includegraphics[width=1.00\textwidth]{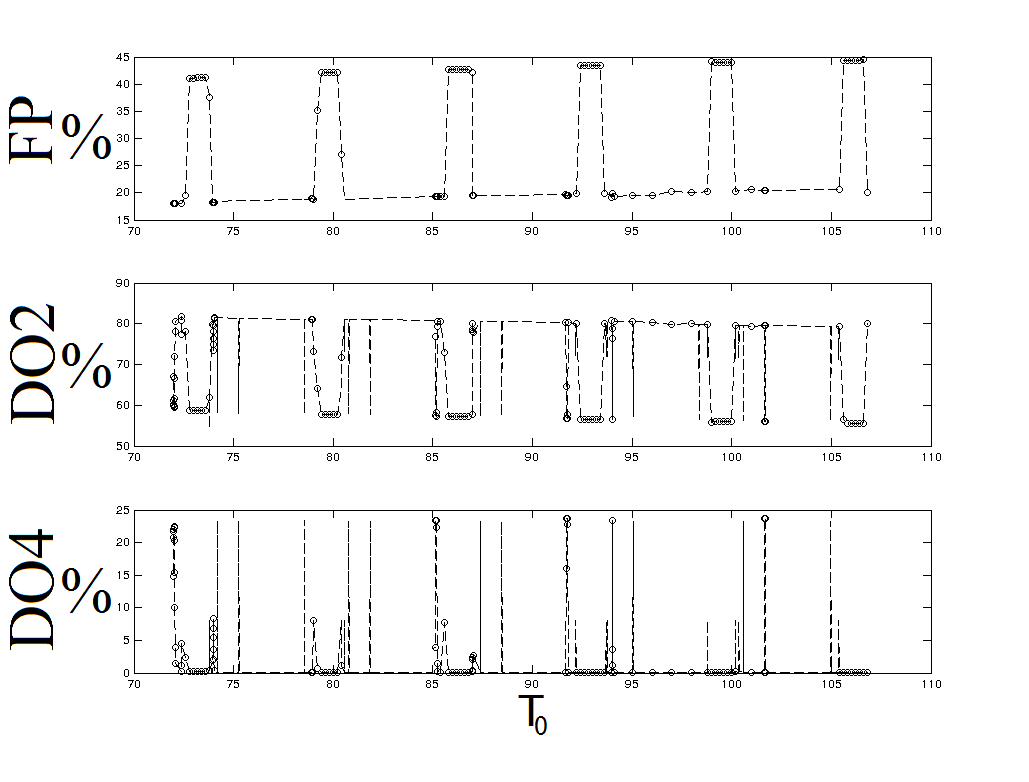}
\vspace{-1.2cm}
\caption{\small Relative areas of the basins of attraction for \eqref{eq:1.2} with
$\alpha = -0.1$, $\beta = 0.545$ and $\gamma(t)$ varying from $ \gamma_i = 0.02$ to $\gamma_f = 0.05$
over a time $T_0$. The circles represent results from numerical simulation, while the dashed lines are
extrapolated on the basis of the observed periodicity --- see the text for details.}
\label{BasinSizes1}
\end{figure}

We interpret the results above as follows. After some transient behaviour, the original sample phase space
has contracted enough and been re-organised so that some regions of phase space become more dense than others;
in particular most of the trajectories end up inside some small, well separated clusters.
The clusters U were converging towards
the attractors PR/NR; however, because of the different evolution of $\gamma(t)$,
they occupy slightly different positions at time $T_1$ as $T_0$ is varied. Since they fall
in a region where the basins of attraction for $\gamma=0.2725$ are sparse and formed of very thin bands,
changing the value of $T_0$ may cause the clusters to cross the boundaries separating the bands belonging
to different basins --- see Figure \ref{DenseRegionsAll}(b).
This results in a large jump in the relative area of the basin of attraction of FP,
as shown in Figure \ref{BasinJumps}; see also Figure \ref{BasinsNotConst}.
According to the results reported in Figure \ref{BasinSizes1},
the clusters U fall inside the basin of attraction of FP for $T_0=93$ and 99,
inside the basin of attraction of DO2 for $T_0=92$, $95$, $96$, $97$ and $98$, and
inside the basin of attraction of DO4 for $T_0=94$.

\vspace{-.4cm}

\vspace{-.4cm}

We expect the jumps to occur for many other values of $T_0$. Figure \ref{BasinSizes1} shows that jumps upward
take place for FP, for instance, at $T_0 \approx 72.6$, 79.2, 85.8, 92.4, 99.0 and 105.6. When a peak appears,
it survives for a very narrow range of $T_0$ (proportional to the width of the bands of the basins of attraction),
after which a jump downward follows.
For this reason, in practice it is difficult to predict exactly the values of $T_0$ where jumps occur,
even though apparently they follow a periodic pattern --- at least in the range of values we have investigated.
Indeed, every value of $T_0$ at which a jump appears is obtained by adding the same quantity $\Delta T_0 \approx 6.6$
to the previous one. For some values of $T_0$ the basin of attraction of the period-4 oscillation
has grown, which implies that the basin of attraction of the new attractor DO4 has
formed exactly where the dense region U has arrived at that time. Also in the case of DO2 and DO4 the jumps
seem to have some periodicity in $T_0$, however with a more intricate structure.
It would be interesting to investigate further such a phenomenon.

Note that, if one applies the method of fast numerical computation described in Section \ref{sec:3}
to compute the relative areas of the basins of attraction in these cases, one needs to take
a much larger value for the time $T_1$ at which the trajectories are approximated to the closest
points of the mesh. More generally, this happens every time (i) an attractor disappears
with increasing dissipation and (ii) the trajectories which were moving towards that attractor
cluster into a small set in a region with tiny bands of the basins of attraction for $\gamma=\gamma_f$.
Then $T_1$ must be large enough for the clusters to have reached
the core surrounding some persisting attractor. When $T_0$ is such that at time $t=2\pi N$
the clusters (which have finite size) touch or even spread across the boundaries of the bands,
time $T_1$ must be very large for the method of fast numerical integration to work, in principle as large
as $T_f$ itself. In this case the method may have no advantage over full time integration of the equations.

\section{Conclusions}

The method of fast numerical computation can be used to estimate the basins of attraction in
systems with damping coefficient $\gamma(t)$ varying linearly between two given values
$\gamma_i$ and $\gamma_f$ over a time $T_0$. First, one fixes a mesh of points in phase space
and for each one computes which basin of attraction it belongs to when the damping coefficient
is constant and equals $\gamma_f$, integrating the equation up to a long time $T_f$.
Then one considers the system with varying dissipation:
one integrates the equations over a time interval $[0,T_1]$ for a large set of initial data
and then approximates each trajectory at time $T_1$ with the closest point of the mesh.
If $T_1$ can be taken much smaller than $T_{f}$, the computational time
can be reduced significantly --- up to a factor 20 in some cases we considered.

The method works well for systems which have at least some regions surrounding the persisting attractors
which are not sensitive to slight perturbations in initial conditions. It becomes more accurate
as the integration time $T_1$ increases, a result of the phase space contracting into regions distant
from the boundaries of the basins of attraction. Thus the results can be improved by setting a
minimum time $T_{\rm min}$, and taking $T_1 =2\pi N \ge \max\{T_0,T_{\rm min}\}$.
This seems to be necessary when $T_0$ is rather small or such that an attractor disappears in a region
where the basins of attraction of the persistent attractors have a multi-band structure.

Using a finer mesh does not necessarily improve the results. However,
smaller basins of attraction benefit from a finer mesh.
Indeed, they intersect only a few points of the mesh and hence, if a coarse mesh is used,
inaccurate predictions are less likely to be cancelled.
The accuracy could be improved further by using a mesh with non-uniform spacing
and better fitting the mesh to the region $C_{T_1}$ to which the sample phase space has contracted at time $T_1$.
This was not implemented here, as it would add extra numerical complexities when rounding a trajectory at
time $T_1$ to the nearest point on the mesh.
Both measures would reduce unnecessary full time integration of initial conditions on the mesh.
Further work could be carried out to remove the use of a mesh and instead have a set of random initial conditions.
Similarly the density of the random initial conditions in particular regions could be chosen so as
to optimise the results.

When the dissipation becomes large, at time $T_1$ the sample region has contracted
into very small clusters. If this happens, either $T_1$ must be very large or a very fine mesh is needed for the method
to produce accurate results. If $T_1$ is not large enough, the step of the mesh has to be
at least comparable with the sizes of the clusters.
For instance, to detect reliably the jumps that may occur in such a case, when an attractor
is destroyed while the damping coefficient increases in time, a too rough approximation of the trajectories
at time $T_1$ would produce a completely wrong description of the dynamics.

Finally, we note that, even though we have considered explicitly in this paper the case where dissipation varies
linearly in time, both for simplicity and for comparison with the literature, we expect our results to extend to more
general situations, where the damping depends on time and either becomes constant after a finite time
or tends asymptotically to a constant value.

\vspace{.4cm}\noindent \textbf{Acknowledgements.}
The integrator ODE113 is one of the built-in ODE solvers in {\sc matlab}.
This research was completed as part of an EPSRC funded PhD.

\end{document}